\documentclass[11pt]{article}
\usepackage{amsfonts, amssymb, amsmath, xcolor, xy, amsthm}
\begin{document}
\newcommand{\QED}{\hfill \blacksquare}
\newcommand{\m}{\textnormal{mult}}
\newcommand{\ra}{\rightarrow}
\newcommand{\la}{\leftarrow}
\renewcommand{\baselinestretch}{1.1}

\theoremstyle{plain}
\newtheorem{thm}{Theorem}[section]
\newtheorem{cor}[thm]{Corollary}
\newtheorem{con}[thm]{Conjecture}
\newtheorem{cla}[thm]{Claim}
\newtheorem{lm}[thm]{Lemma}
\newtheorem{prop}[thm]{Proposition}
\newtheorem{example}[thm]{Example}

\theoremstyle{definition}
\newtheorem{dfn}[thm]{Definition}
\newtheorem{alg}[thm]{Algorithm}
\newtheorem{rem}[thm]{Remark}

\title{\bf  An Analogue of the Gallai-Edmonds Structure Theorem for Nonzero Roots of the Matching Polynomial}
\author{
Cheng Yeaw Ku
\thanks{ Department of Mathematics, Caltech, Pasadena, CA 91125, USA. E-mail: cyk@caltech.edu.} \and William Chen \thanks{ MSC 176, Caltech, Pasadena, CA 91125, USA. E-mail:
chenw@caltech.edu} }

\maketitle

\begin{abstract}\noindent
Godsil observed the simple fact that the multiplicity of 0 as a root of the matching polynomial of a graph coincides with the classical notion of deficiency. From this fact he asked to what extent classical results in matching theory generalize, replacing ``deficiency" with multiplicity of $\theta$ as a root of the matching polynomial. We prove an analogue of the Stability Lemma for any given root, which describes how the matching structure of a graph changes upon deletion of a single vertex. An analogue of Gallai's Lemma follows. Together these two results imply an analogue of the Gallai-Edmonds Structure Theorem. Consequently, the matching polynomial of a vertex transitive graph has simple roots.
\end{abstract}

\section{Introduction}
A {\em matching} of a graph $G$ is a set of pairwise non-adjacent edges of $G$. Classical matching theory is mostly concerned with the maximum size $\nu(G)$ of a matching in $G$, known as the {\it matching number}. Another important quantity is the number of vertices $\mbox{def}(G)$ missed by a maximum matching, known as the {\it deficiency}. They are related by the formula $\mbox{def}(G)=|V(G)|-2\nu(G)$.

Recall that for a graph $G$ on $n$ vertices, the {\it matching
polynomial} $\mu(G,x)$ of $G$ is given by
\[ \mu(G,x)=\sum_{k \ge 0} (-1)^{k}p(G,k)x^{n-2k}, \]
where $p(G,k)$ is the number of matchings with $k$ edges in $G$. Let $\m(\theta, G)$ denote the multiplicity
of $\theta$ as a root of $\mu(G,x)$. Godsil observed that $\m(0,G)=\mbox{def}(G)$ and obtained several results in \cite{G} generalizing classical results involving $\mbox{def}(G)$ to analogous $\m(\theta,G)$ versions. The present paper is concerned with such a generalization for the celebrated Gallai-Edmonds Structure Theorem.

The following definition introduced by Godsil in \cite{G} (who adapted it from \cite{N}) is useful in stating the theorem. It assigns to each vertex a ``sign''---minus, zero, or plus---based on how the multiplicity of a root of the matching polynomial changes when that vertex is deleted. This definition is fundamental to the work here.
\begin{dfn}\label{def}
Let $\theta$ be a root of $\mu(G,x)$. For any vertex $u\in V(G)$,
\begin{itemize}
\item $u$ is $\theta$-{\em essential} if $\m(\theta, G
\setminus u)<\m(\theta, G)$,

\item $u$ is $\theta$-{\em neutral} if $\m(\theta, G
\setminus u)=\m(\theta, G)$,

\item $u$ is $\theta$-{\em positive} if $\m(\theta, G
\setminus u)>\m(\theta, G)$.
\end{itemize}
\end{dfn}
\begin{rem}A vertex is $0$-essential if and only if it is missed by
some maximum matching of $G$. There are no $0$-neutral vertices.
\end{rem}
\begin{rem}If $\m(\theta, G)=0$ then there are no $\theta$-essential
vertices since the multiplicity of a root cannot be negative.
Nevertheless, it still makes sense to talk about $\theta$-neutral
and $\theta$-positive vertices when $\m(\theta, G)=0$.
\end{rem}
We will often omit the $\theta$- prefix from these terms if it is clear from context.

A further classification of vertices plays an important role in describing the Gallai-Edmonds Structure Theorem:
\begin{dfn}Let $\theta$ be a root of $\mu(G,x)$. For any vertex $u\in V(G)$, $u$ is $\theta$-{\em special} if it is
not $\theta$-essential but has a neighbor that is $\theta$-essential.
\end{dfn}

The Gallai-Edmonds Structure Theorem describes a certain canonical decomposition of $V(G)$. Its statement
essentially consists of two lemmas, the Stability Lemma and Gallai's Lemma. For more
information, see \cite[Section 3.2]{LP}. The main results of the present paper are the following analogues for the Stability Lemma and Gallai's Lemma for any root of the matching polynomial. Their classical counterparts are simply the case $\theta=0$.
\begin{thm}[Stability Lemma]\label{main}
\noindent Let $G$ be a graph with $\theta$ a root of
$\mu(G,x)$, $u$ a $\theta$-special vertex in
$G$ and $v$ a vertex of $G$ different from $u$. Then
\begin{itemize}
\item $v$ is $\theta$-essential in $G$ if and only if $v$ is $\theta$-essential in $G\setminus u$,

\item $v$ is $\theta$-neutral in $G$ if and only if $v$ is $\theta$-neutral in $G\setminus u$,

\item $v$ is $\theta$-positive in $G$ if and only if $v$ is $\theta$-positive in $G\setminus u$.
\end{itemize}
\end{thm}
\begin{rem}
This result is slightly different from the classical Stability Lemma because it includes neutral vertices. Recall that there are no $0$-neutral vertices.
\end{rem}

It follows from Theorem \ref{main} that after deleting the special vertices, the essential vertices remain essential. Furthermore, they are not joined to the other non-essential vertices. Therefore it is interesting to study graphs whose vertices are all $\theta$-essential. These graphs are called {\it $\theta$-primitive}, and generalize the factor-critical graphs. Gallai's Lemma is a fundamental result about the structure of these graphs when $\theta=0$. We prove this for any root $\theta$.
\begin{thm}[Gallai's Lemma]\label{main2}
\noindent If every vertex of $G$ is $\theta$-essential, then $\m(\theta,G)=1$.
\end{thm}

The following corollary is immediate since every vertex of a vertex transitive graph is $\theta$-essential, see \cite{G}.
\begin{cor}
The matching polynomial of a vertex transitive graph has simple roots.
\end{cor}

This answers a question of Godsil in \cite[Problem 6.1]{G1} and disproves a conjecture of Mohar \cite{M}: for every integer $r$ there exists a vertex transitive graph $G$ whose matching polynomial has a root of multiplicity at least $r$.

\section{Basic properties}
In this section, we collect some basic identities and properties of the matching polynomial proved in \cite{G0} and \cite{G}. If $u\in V(G)$, then $G \setminus u$ is the graph obtained from $G$ by deleting vertex $u$ and the edges of $G$ incident to $u$. We also denote the graph $(G \setminus u) \setminus v$ by $G \setminus uv$. If $e\in E(G)$, the graph $G-e$ is the graph obtained from $G$ by deleting the edge $e$. If $f\not\in E(G)$ is a pair of distinct vertices, then $G+f$ is the graph obtained by adding $f$ as an edge to $G$.

The matching polynomial satisfies the following basic identities.
\begin{prop}
Let $G$ and $H$ be graphs, with matching polynomials $\mu(G,x)$ and $\mu(H,x)$, respectively. Then
\begin{itemize}
\item[\textnormal{(a)}] $\mu(G \cup H,x) = \mu(G,x)\mu(H,x)$,
\item[\textnormal{(b)}] $\mu(G,x)=\mu(G - e, x)-\mu(G
\setminus uv,x)$ where $e=\{u,v\}$ is an edge of $G$,
\item[\textnormal{(c)}] $\mu(G,x) = x\mu(G \setminus u,x)-\sum_{v
\sim u} \mu(G \setminus uv,x)$ for any vertex $u$ of $G$.
\end{itemize}
\label{identity}
\end{prop}
Proposition \ref{identity}(a) says that the matching polynomial can be considered separately for each connected component of a disconnected graph. We will use Proposition \ref{identity}(b) frequently, and it is especially applicable to the results of Section 3.

Another useful result due to Godsil guarantees the existence of a $\theta$-essential vertex in a graph whose matching polynomial has $\theta$ as a root. This implies that every vertex of a vertex transitive graph is $\theta$-essential for any root $\theta$.
\begin{lm}\label{essential}
Any $G$ with $\m(\theta, G)>0$ must have at
least one $\theta$-essential vertex.
\end{lm}

Godsil \cite{G0} showed that the roots of $G\setminus u$ interlace those of $G$. This puts a limitation on how much the multiplicity of a given root of the matching polynomial can change upon deleting a vertex. In this paper, as in Godsil's, we often refer to this phenomenon as ``interlacing.''
\begin{prop}
Let $G$ be a graph, $u \in V(G)$ a vertex of $G$. Then $\m(\theta,G\setminus u)$ differs from $\m(\theta,G)$ by at most one.
\end{prop}
The notions of $\theta$-essential, neutral, and positive introduced in Definition \ref{def} should be viewed under this useful proposition.

The next result says that a $\theta$-special vertex must be $\theta$-positive and has significant consequences for the Gallai-Edmonds decomposition.
\begin{lm}\label{no-neutral-essential}
A $\theta$-neutral vertex cannot be joined to a $\theta$-essential
vertex.
\end{lm}

If $P$ is a path in $G$, then $G \setminus P$ denotes the graph
obtained from $G$ be deleting the vertices of $P$ and all the edges
incident to these vertices.

One of its important corollaries implies that if $u$ and $v$ are adjacent vertices of $G$, then $u$ cannot be $\theta$-essential in $G\setminus v$ if $v$ is $\theta$-essential in $G$. More generally,
\begin{cor}\label{interlacing-path}
For any root $\theta$ of $\mu(G,x)$ and a path $P$ in $G$,
\[ \m(\theta, G \setminus P) \ge \m(\theta, G)-1. \]
\end{cor}
As in Godsil's paper \cite{G}, we say that $P$ is {\it $\theta$-essential} if $\m(\theta, G \setminus P)
= \m(\theta, G)-1$. This connection is further motivated by the following.
\begin{lm}\label{no-path}
The end vertices of a $\theta$-essential path are themselves $\theta$-essential.
\end{lm}

Using these tools, Godsil proved a result very similar to the Stability Lemma.
\begin{prop}[Theorem 4.2, \cite{G}]\label{delete-positive}
Let $\theta$ be a root of $\mu(G,x)$
and let $u$ be a $\theta$-positive vertex in $G$. Then
\begin{itemize}
\item[\textnormal{(a)}] if $v$ is $\theta$-essential in $G$ then it
is $\theta$-essential in $G \setminus u$,
\item[\textnormal{(b)}] if $v$ is $\theta$-positive in $G$ then it
is $\theta$-essential or $\theta$-positive in $G \setminus u$,
\item[\textnormal{(c)}] if $u$ is $\theta$-neutral in $G$ then it is
$\theta$-essential or $\theta$-neutral in $G \setminus u$.
\end{itemize}
\end{prop}
\begin{rem}
The Stability Lemma says that the sign of a vertex does not change upon deleting a special vertex. Proposition \ref{delete-positive} investigates how the sign changes when deleting a positive vertex.
\end{rem}

It is not difficult to formulate similar results for neutral vertices using the same techniques from Godsil's proof in \cite{G} of Proposition \ref{delete-positive}. In many ways, positive and neutral vertices behave similarly. Since the proof is almost identical to that of Godsil, we omit it here.

\begin{prop}\label{delete-neutral}
Let $\theta$ be a root of $\mu(G,x)$
and let $u$ be a $\theta$-neutral vertex in $G$. Then

\noindent (a) if $v$ is $\theta$-essential in $G$ then it is
$\theta$-essential in $G \setminus u$,

\noindent (b) if $v$ is $\theta$-positive in $G$ then it is
$\theta$-positive or $\theta$-neutral in $G \setminus u$,

\noindent (c) if $v$ is $\theta$-neutral in $G$ then it is
$\theta$-neutral or $\theta$-positive in $G \setminus u$.
\end{prop}

The result for essential vertices follows easily from the previous two.

\begin{prop}\label{delete-essential}
Let $\theta$ be a root of $\mu(G,x)$
and let $u$ be a $\theta$-essential vertex in $G$.

\noindent (a) if $v$ is $\theta$-positive in $G$ then it is
$\theta$-positive in $G \setminus u$,

\noindent (b) if $v$ is $\theta$-neutral in $G$ then it is
$\theta$-neutral in $G \setminus u$.

\noindent In particular, if $v$ is $\theta$-essential in $G
\setminus u$ where $u$ is $\theta$-essential in $G$, then $v$ is
$\theta$-essential in $G$.
\end{prop}
\begin{proof} Suppose $v$ is $\theta$-positive in $G$. Then, by Proposition
\ref{delete-positive}, $\m(\theta, G \setminus uv)=\m(\theta, G
\setminus vu)=\m(\theta, G)$, so $v$ is $\theta$-positive in $G \setminus u$.
Now, suppose $v$ is $\theta$-neutral in $G$. By Proposition
\ref{delete-neutral}, $\m(\theta, G \setminus uv) = \m(\theta, G
\setminus vu) = \m(\theta, G)-1$ so that $v$ is neutral in $G \setminus u$.
\end{proof}
The proof of Proposition \ref{delete-essential} is based on the trivial observation that the order in which vertices are deleted is immaterial, a technique that is used extensively in this paper.
\begin{rem}
Propositions \ref{delete-positive}, \ref{delete-neutral}, and \ref{delete-essential} are best possible in the sense that they place the most severe restrictions on the sign of the vertices of $G\setminus u$ in each case. That is, only the possibilities explicitly excluded do not occur.
\end{rem}

\section{Edge manipulations}

Let $G^*$ be the graph obtained by adding an edge to $G$, say $f=\{u,v\}$. Since $G^*\setminus u=G\setminus u$ and $G^*\setminus v=G\setminus v$, the signs of $u$ and $v$ must be the same relative to each other. The actual signs are determined by the multiplicity of $G^*$. The same argument works when deleting edges.

First we consider what happens to the multiplicity of $\theta$ upon adding an edge.

\begin{lm}\label{add-essential-positive}
Let $u$ be a $\theta$-positive vertex in $G$. Then for any nonadjacent vertex $v\ne u$,
$$\m(\theta, G+f) = \m(\theta, G),$$
where $f=\{u,v\} \not \in E(G)$.
Therefore $u$ is $\theta$-positive in $G+f$ and $v$ has the same sign in $G+f$ as it did in $G$.
\end{lm}
\begin{proof} Let $k=\m(\theta, G)$ and $G^*=G+f$. Recall the statement of Proposition \ref{identity}(b), which in this case states that
\begin{equation}
\mu(G^*,x)=\mu(G,x)-\mu(G^*\setminus uv,x).
\label{eq1}
\end{equation}
Since $u$ is positive, $\m(\theta,G^*\setminus uv)=\m(\theta,G\setminus uv)\geq k$, and (\ref{eq1}) gives that $\m(\theta,G^*)\geq k$.

If $v$ is essential in $G$, $\m(\theta,G^*\setminus v)=\m(\theta,G\setminus v)=k-1$, so by interlacing $\m(\theta, G^*)\leq k$.

If $v$ is neutral in $G$, $\m(\theta,G^*\setminus v)=k$ so by interlacing $\m(\theta, G^*)\leq k+1$. If $\m(\theta, G^*)=k+1$ then $u$ is neutral and $v$ is essential in $G^*$, contradicting Lemma \ref{no-neutral-essential}. It follows that $\m(\theta, G^*)\leq k$.

If $v$ is positive in $G$, then by Proposition \ref{delete-positive}, either $\m(\theta,G^*\setminus uv)=k+2$ or $\m(\theta,G^*\setminus uv)=k$. In the first case, $\m(\theta, G^*) \le k$ by (\ref{eq1}) and we are done. In the second case, $v$ is essential in $G^*\setminus u$. But this is impossible, because if $u, v$ are both neutral in $G^*$ then this contradicts Proposition \ref{delete-neutral}, and if $u,v$ are both essential in  $G^*$ then this contradicts Corollary \ref{interlacing-path}.
\end{proof}

\begin{lm}\label{add-neutral-essential}
Let $u$ be a $\theta$-neutral vertex and $v \ne u$ be a nonadjacent $\theta$-essential
vertex in $G$. Then $\m(\theta, G+f)= \m(\theta, G)-1$, where $f=\{u,v\} \not \in E(G)$. Therefore $u$ is $\theta$-positive and $v$
is $\theta$-neutral in $G+f$.
\end{lm}
\begin{proof} Let $k=\m(\theta, G)$ and $G^*=G+f$. By Proposition \ref{delete-neutral}, $\m(\theta, G^*\setminus uv)=\m(\theta, G \setminus uv)=k-1$. By (\ref{eq1}) and interlacing, it follows that $\m(\theta, G^*)=k-1$.
\end{proof}

For the other cases the situation is not as clean. Of those cases, the following lemma will be useful for our purposes, although similar results can be proven for other sign combinations.
\begin{lm}\label{add-essential-essential-I}
Let $u, v$ be nonadjacent $\theta$-essential vertices in $G$ such that $\m(\theta, G \setminus uv)
\ge \m(\theta, G)-1$. Then, either
\begin{itemize}
\item[$\bullet$] $\m(\theta, G+f)=\m(\theta, G)-1$ and both $u$ and $v$ are $\theta$-neutral in $G+f$, or
\item[$\bullet$] $\m(\theta, G+f)=\m(\theta, G)$ and both $u$ and $v$ are $\theta$-essential in $G+f$.
\end{itemize}
\end{lm}
\begin{proof} Let $k=\m(\theta, G)$ and $G^*=G+f$. By (\ref{eq1}), $\m(\theta, G^*) \ge k-1$
using the assumption that $\m(\theta, G \setminus uv) \ge k-1$. Since
$\m(\theta, G^* \setminus u)=\m(\theta, G \setminus u)=k-1$, by interlacing it follows that
$\m(\theta, G^*) \le k$.
\end{proof}

Now we consider what happens to
the multiplicity of $\theta$ when we
delete an edge $e=\{u,v\}$ from $G$.

\begin{lm}\label{delete-edge}
Let $u$ be a $\theta$-special vertex in $G$, adjacent to a
$\theta$-essential vertex $v$. Let $e=\{u,v\} \in E(G)$. Then $\m(\theta, G-e) = \m(\theta, G)$, therefore $u$ remains
$\theta$-positive and $v$ remains $\theta$-essential in $G-e$.
\end{lm}
\begin{proof}
Let $k=\m(\theta, G)$ and $G'=G-e$. Notice that $\m(\theta,G'\setminus u)=\m(\theta,G\setminus u)=k+1$ and $\m(\theta,G'\setminus v)=\m(\theta,G\setminus v)=k-1$. By interlacing it follows that $\m(\theta,G)=k$.
\end{proof}

\begin{lm}\label{delete-neutral-positive}
Let $u$ be a $\theta$-positive vertex in $G$, adjacent to a $\theta$-neutral
vertex $v$. Let $e=\{u,v\} \in E(G)$. Then, either
\begin{itemize}
\item[$\bullet$] $\m(\theta, G-e) = \m(\theta, G)+1$, $u$ is $\theta$-neutral and $v$ is $\theta$-essential in $G-e$, or
\item[$\bullet$] $\m(\theta, G-e) = \m(\theta, G)$, $u$ is $\theta$-positive and $v$ is $\theta$-neutral in $G-e$.
\end{itemize}
\end{lm}

\begin{proof}Let $k=\m(\theta, G)$ and $G'=G-e$. By Proposition
\ref{delete-positive}, $\m(\theta, G \setminus uv) \ge k$.
Applying \ref{identity}(b), we have
$\m(\theta, G') \ge k$. As $\m(\theta, G' \setminus v) = \m(\theta,
G \setminus v) = k$, it follows that $\m(\theta, G') \le k+1$ by
interlacing.
\end{proof}

\section{Three lemmas}

In this section, we study the effect of deleting an edge
incident to a $\theta$-special vertex. This will yield three lemmas used in the proof of Theorem \ref{main} by induction. We first consider the case
when a $\theta$-special vertex has two $\theta$-essential neighbors.

\begin{lm}\label{bac-not-essential}
Let $u$ be a $\theta$-special vertex in $G$ adjacent to two
$\theta$-essential vertices $v$ and $w$ in $G$, and let $e=\{u,v\} \in E(G)$. Suppose that the path $vuw$ is not $\theta$-essential in $G$. Then $u$ is $\theta$-special in
$G-e$, $w$ is $\theta$-essential in $G-e$ and
$\m(\theta, G-e)=\m(\theta, G)$.
\end{lm}

\begin{proof}Let $G'=G-e$ and $k=\m(\theta, G)$. By Lemma \ref{delete-edge}, it follows that
$\m(\theta, G')=k$, $u$ is positive in $G'$ and $v$ is
essential in $G'$, so it is enough to show that $w$ remains an
essential neighbor of $u$ in $G'$.

Notice by Proposition \ref{delete-essential} $u$ is positive in $G\setminus w$. Also, $v$ cannot be essential in $G\setminus w$, otherwise by Proposition \ref{delete-positive} the path $vuw$ is essential in $G$. So $v$ is either neutral or positive in $G\setminus w$.

If $v$ is neutral in $G\setminus w$, then by Lemma \ref{delete-neutral-positive} it follows that either $\m(\theta,G'\setminus w)=k$ or $\m(\theta,G'\setminus w)=k-1$. In the latter case we are done, so we show that the first case is not possible. In that case, $u$ is neutral and $v$ is essential in $G'\setminus w$, so by Proposition \ref{delete-neutral} $\m(\theta,G\setminus vuw)=\m(\theta,G'\setminus wuv)=k-1$, contradicting the assumption that the path $vuw$ is not essential in $G$.

If $v$ is positive in $G\setminus w$, then $u$ must be positive in $G\setminus wv$, otherwise by Proposition \ref{delete-positive} $u$ is essential in $G\setminus wv$ so $vuw$ is an essential path in $G$. Therefore, $\m(\theta,G\setminus vuw)=k+1$. Now consider the sign of $w$ in $G'$. The vertex $w$ cannot be neutral in $G'$, otherwise $\m(\theta,G\setminus wv)=\m(\theta, G'\setminus wv)=k-1$ by Proposition \ref{delete-neutral} so $\m(\theta,G\setminus vuw)\ne k+1$ by interlacing. If $w$ is essential in $G'$ we are done, so we may assume $w$ is positive in $G'$.

Since $\m(\theta, G' \setminus
wu)=\m(\theta, G \setminus uw)=k$, $u$ is
essential in $G' \setminus w$. By Proposition \ref{delete-positive}, $v$ is also
essential in $G' \setminus w$. Since $\m(\theta, G' \setminus w)=k+1$, applying Lemma \ref{add-essential-essential-I} to
$G' \setminus w$ yields $\m(\theta, G \setminus w)=\m(\theta, (G' \setminus w)+e) \ge k$, contradicting that $w$ is essential in $G$.
\end{proof}

Next, we consider the situation in which a $\theta$-special vertex
$u$ has a $\theta$-essential neighbor $v$ and a $\theta$-neutral
neighbor $w$. It turns out that $u$ is still $\theta$-special after deleting the edge $\{u,w\}$.

\begin{lm}\label{special-neutral}
Let $u$ be a $\theta$-special vertex and $v$ be a $\theta$-essential
neighbor of $u$ in $G$. Suppose $w$ is a $\theta$-neutral neighbor
of $u$ in $G$, $e=\{u,w\} \in E(G)$. Then $u$ is $\theta$-special in $G-e$, $v$
is $\theta$-essential in $G-e$ and $\m(\theta, G-e)=\m(\theta, G)$.
\end{lm}

\begin{proof} Let $G'=G-e$ and $k=\m(\theta, G)$. By Lemma
\ref{delete-neutral-positive}, either $\m(\theta, G')=k+1$ or $\m(\theta, G')=k$.

If $\m(\theta, G')=k+1$, then $u$ is
neutral and $w$ is essential in $G'$. Since $\m(\theta, G' \setminus u) = \m(\theta, G \setminus u) =
k+1$ and $\m(\theta, G' \setminus uv)=\m(\theta, G \setminus uv) =
k$, $v$ must be essential in $G' \setminus u$. As $u$ is
neutral in $G'$, by Proposition \ref{delete-neutral}, $v$
must be essential in $G'$, contradicting Lemma
\ref{no-neutral-essential}.

If $\m(\theta, G')=k$, then $u$ is positive and $w$ is neutral in $G'$. By Proposition \ref{delete-neutral} $\m(\theta, G' \setminus wv) = \m(\theta, G \setminus vw)=k-1$. So $v$ is
essential in $G' \setminus w$. As $w$ is neutral
in $G'$, by Proposition \ref{delete-neutral} again, $v$ is
essential in $G'$. So $u$ is special in $G'$ since
it is positive in $G'$ and is joined to $v$ in $G'$.
\end{proof}

A similar result holds when $u$ is adjacent to a $\theta$-positive vertex.

\begin{lm}\label{special-positive}
Let $u$ be a $\theta$-special vertex in $G$ and $v$ a
$\theta$-essential neighbor of $u$ in $G$. Suppose $w$ is a $\theta$-positive neighbor of $u$
in $G$, $e=\{u,w\}$. Then $u$ is $\theta$-special in $G-e$, $v$ is
$\theta$-essential in $G-e$ and $\m(\theta, G-e)=\m(\theta,G)$.
\end{lm}

\begin{proof}
Let $G'=G-e$ and $k=\m(\theta, G)$.

If $u$ were neutral in $G'$, then $\m(\theta,
G')=k+1$. By Lemma \ref{no-neutral-essential}, $v$ cannot be
essential in $G'$. So, by Proposition \ref{delete-neutral},
we have $\m(\theta, G' \setminus uv)
\ge k+1$, contradicting that $\m(\theta, G' \setminus uv)=\m(\theta, G \setminus uv)=k$. So $u$ cannot be neutral in $G'$. If $u$ were essential in $G'$,
then $\m(\theta,
G')=k+2$. But $\m(\theta,G'\setminus uv)=\m(\theta,G\setminus uv)=k$, contradicting Corollary \ref{interlacing-path}.

Therefore $u$ is positive in $G'$, and $\m(\theta, G')=k$. Using Lemma \ref{identity} (b),
\begin{eqnarray}\mu(G
\setminus v, x)=\mu(G' \setminus v, x)-\mu(G \setminus vuw, x) \label{e1}.
\end{eqnarray}
If $v$ is not essential in $G'$ then $\m(\theta, G' \setminus v) \ge k$,
so by Lemma \ref{no-path}, the multiplicity of $\theta$ on the right hand side of (\ref{e1})
is always at least $k$, contradicting the fact that $\m(\theta, G
\setminus v)=k-1$ on the left hand side. Therefore, $v$ is
essential in $G'$ and so $u$ is special in $G'$.
\end{proof}
\section{The Gallai-Edmonds Structure Theorem}

We are now ready to prove Theorem \ref{main}. In view of Proposition
\ref{delete-positive}, it remains to show that for any $\theta$-special vertex $u$, $v$ is
$\theta$-essential in $G \setminus u$ only if $v$ is $\theta$-essential
in $G$.

It is easy to show that $v$ cannot be $\theta$-neutral in $G$.
\begin{prop}\label{not-neutral}
Suppose $u$ is $\theta$-special in $G$ and $v$ is $\theta$-essential in $G \setminus u$. Then $v$ cannot be $\theta$-neutral in $G$.
\end{prop}
\begin{proof}
Suppose $v$ is neutral in $G$ and $k=\m(\theta,G)$. Let $w$ be an essential vertex adjacent to $u$ in $G$. Since $\m(\theta,G\setminus uv)=k$, $u$ is neutral in $G\setminus v$. But $w$ is essential in $G\setminus v$, contradicting Lemma \ref{no-neutral-essential}.
\end{proof}

\begin{proof}[Proof of Theorem \ref{main}]
The proof is by induction on the degree of $u$. Let $w_1$ be an essential vertex witnessing that $u$ is special, and let $e=\{u,w_1\} \in E(G)$. We may also assume that $\theta \not = 0$.

\noindent{\bf Base Case:} Notice that $w_1$ cannot be the only neighbor of $u$, otherwise by Lemma \ref{delete-edge}, $u$ is positive in $G-e$, a contradiction since $u$ is isolated in $G-e$. Suppose $\deg(u)=2$. Let $w_2$ be the second neighbor of $u$. Then by Lemmas \ref{bac-not-essential}, \ref{special-neutral}, and \ref{special-positive}, it follows that $w_2$ is essential and the path $w_1uw_2$ is essential in $G$, otherwise by deleting an edge $u$ would be a special vertex with only one neighbor.

For a contradiction, we now assume, in view of Proposition \ref{not-neutral}, that $v$ is positive in
$G$. We first prove the following claims. Let $G'=G - e$. Note that $\m(\theta, G')=k$ by Lemma \ref{delete-edge}.

\noindent {\bf Claim 1.}~~$v$ is positive in $G'$.

Let $G^{*}=G+f$ where $f=\{v,u\} \not \in E(G)$. Since $G^{*} \setminus u =
G \setminus u$ and $G^{*} \setminus v = G \setminus v$, $u$ and $v$
must be both essential or both positive or both
neutral in $G^{*}$. Suppose $u$ and $v$ are both
essential or both neutral in $G^{*}$. Then
$\m(\theta, G^{*}) \ge k+1$ since $\m(\theta, G^{*} \setminus u) =
\m(\theta, G \setminus u) =k+1$. By the interlacing property of a
path, $\m(\theta, G^{*} \setminus w_1uw_2) \ge k$, contradicting the
fact that $\m(\theta, G^{*} \setminus w_1uw_2) = \m(\theta, G \setminus
w_1uw_2)=k-1$. Therefore, $u$ and $v$ are both positive in
$G^{*}$ and $\m(\theta, G^{*})=k$.

If $v$ is neutral in $G'$ then $\m(\theta, G' \setminus vw_1)
= k-1$ by Proposition \ref{delete-neutral}. But $\m(\theta, G'
\setminus vw_1)=\m(\theta, G \setminus vw_1)=k$ by Proposition
\ref{delete-positive} since $v$ is positive and $w_1$ is
essential in $G$. So $v$ is either positive or
essential in $G'$.

Suppose $v$ is essential in $G'$. By Proposition
\ref{delete-positive}, $w_2$ is essential in $G \setminus u =
G' \setminus u$. Let $H$ denote the graph which is the union of $G'
\setminus u$ and the isolated vertex $u$. By Lemma
\ref{add-neutral-essential}, $w_2$ is neutral in $G'=H+\{u,w_{2}\}$, $u$ is
positive in $G'$ and $\m(\theta, G')=k$. Now, deleting $v$
first from $G'$ followed by deleting $w_2$ and $u$, we deduce that
$\m(\theta, G' \setminus vw_2u)=k-1$ using Proposition
\ref{delete-essential} and Proposition \ref{delete-neutral}. Since
$G^{*} \setminus vuw_2 = G' \setminus vuw_2$, we deduce that $vuw_2$ is
essential in $G^{*}$, whence $v$ is essential in
$G^{*}$ by Lemma \ref{no-path}. This contradicts the conclusion of
the first paragraph following Claim 1.

Hence, $v$ is positive in $G'$, proving Claim 1.

\noindent {\bf Claim 2.}~~$\m(\theta, G' \setminus vw_2)=k$.

Since $\m(\theta, G' \setminus vu)=\m(\theta, G \setminus uv)=k$
(recall that $v$ is essential in $G \setminus u$ by our
hypothesis), it follows immediately from Claim 1 that $u$ is
essential in $G' \setminus v$. Then, being adjacent to $u$,
$w_{2}$ is either essential or positive in $G'
\setminus v$ (Lemma \ref{no-neutral-essential}). If $w_{2}$ is
positive in $G' \setminus v$, then $u$ remains
essential in $(G' \setminus v)-e'$ where $e'=uw_2 \in E(G'
\setminus v)$ by Lemma \ref{delete-edge}. However, as an isolated
vertex in $(G' \setminus v) - e'$, $u$ has to be neutral
(since $\theta \not =0$) in $(G' \setminus v) - e'$, contradicting
the preceding sentence. So $w_{2}$ is essential in $G'
\setminus v$ and $\m(\theta, G' \setminus vw_2)=k$, thus proving Claim
2.

Finally, recall that $w_1$ is essential in $G'$ (Lemma
\ref{delete-edge}) and $v$ is positive in $G'$ (Claim 1).
Also, by Lemma \ref{add-neutral-essential}, $w_2$ is
neutral in $G'=H+\{u,w_{2}\}$. By Claim 2, $v$ is neutral in $G' \setminus w_{2}$. Clearly, as $\theta \not = 0$, the isolated vertex
$u$ is neutral in $G' \setminus vw_2$. Subsequently, using
Proposition \ref{delete-neutral}, by deleting $w_{2}$ from $G'$ followed by deleting $v$ from $G'\setminus w_{2}$, we deduce that $\m(\theta, G \setminus
vw_2uw_1) = \m(\theta, (G' \setminus w_2v) \setminus uw_1)=k-1$. On the
other hand, by interlacing (Corollary \ref{interlacing-path}),
$\m(\theta, G \setminus vw_1uw_2)= \m(\theta, (G \setminus v) \setminus
w_1uw_2) \ge \m(\theta, G \setminus v)-1 = k$ since we assume $v$ is positive in $G$, contradicting the preceding sentence. This establishes the theorem when $\deg(u)=2$.

We may now assume that $\deg(u) \ge 3$.

\noindent{\bf Inductive Step:} Let $w_1$ be an essential neighbor of $u$ witnessing that $u$ is special. If $w_2\neq w_1$ is adjacent to $u$ and $w_1uw_2$ is not a essential path, let $e_2=\{u,w_2\}$ and $G_2=G- e_2$. By the three Lemmas \ref{bac-not-essential}, \ref{special-neutral}, and \ref{special-positive}, $u$ is still special in $G_{2}$ and $\m(\theta, G_2)=\m(\theta, G)$. By the induction hypothesis, $v$ is essential in $\m(\theta, G_{2})$, so $\m(\theta, G_2\setminus v)=\m(\theta,G)-1$. Since $u$ is still positive in $G_2\setminus v$, by Lemma \ref{add-essential-positive} $\m(\theta, G\setminus v)=\m(\theta, (G_2\setminus v)+e_2)=\m(\theta,G)-1$, so $v$ is essential in $G$.

If for every vertex $w\neq w_1$ adjacent to $u$ the path $w_1uw$ is an essential path, let $w_2,w_3$ be two such vertices. Let $e_3=\{u,w_3\}$ and $G_3=G- e_3$. By Lemma \ref{delete-edge}, $\m(\theta,G_3)=\m(\theta,G)$ and $u$ is positive in $G_3$. Since $w_1uw_2$ is still an essential path in $G_3$, $w_1,w_2$ are essential in $G_{3}$ (Lemma \ref{no-path}), so $u$ is special in $G_{3}$. Now the proof follows as before: by the induction hypothesis, $v$ is essential in $G_3$, so $\m(\theta, G_3\setminus v)=\m(\theta,G)-1$. Since $u$ is still positive in $G_3\setminus v$, by Lemma \ref{add-essential-positive} $\m(\theta, G\setminus v)=\m(\theta, (G_2\setminus v)+e_3)=\m(\theta,G)-1$, so $v$ is essential in $G$.
\end{proof}

With the Stability Lemma in hand, we can state a weak version of the Gallai-Edmonds Structure Theorem. Denote by $A(G)$ the set of all $\theta$-special vertices of $G$ for some root $\theta$ of $\mu(G,x)$. Deleting the $\theta$-special vertices one by one, the Stability Lemma says that the $\theta$-essential vertices of $G\setminus A(G)$ form $\theta$-primitive components and, by Lemma \ref{essential}, the non-essential vertices form components not having $\theta$ as a root. Let $D(G)$ be the set of $\theta$-essential vertices of $G$ and $C(G)=V(G)\setminus (A(G)\cup D(G))$. The partition of $V(G)$ into $A(G)$, $C(G)$, and $D(G)$ is called the {\it Gallai-Edmonds decomposition}. It will be useful to keep the above in mind for the proof of Theorem \ref{main2}, which states that if every vertex of a graph $G$ is $\theta$-essential, then $\m(\theta,G)=1$.

\begin{proof}[Proof of Theorem \ref{main2}]
Set $k=\m(\theta,G)$. Assume for a contradiction that $\m(\theta, G \setminus v)=k-1>0$. Using the notation above, let $D=D(G\setminus v)$,  $A=A(G\setminus v)$, and $C=C(G\setminus v)$. Let $w\in A$. Starting from $w$, we now delete the vertices of $A$ from $G$ one by one. The multiplicity of $\theta$ in $G\setminus A$ is at most $k+|A|-2$, since $w$ is essential in $G$ and by interlacing deleting the other $|A|-1$ vertices increases the multiplicity by at most $|A|-1$.

Abusing notation, let $D$ be the subgraph of $G$ induced by $D$. Since $\m(\theta,G\setminus v)=k-1$, it follows from the Stability Lemma applied to $G \setminus v$ that $\m(\theta,D)=k-1+|A|$.  By Corollary \ref{interlacing-path}, $v$ is not adjacent to any vertices of $D$. Since $D$ are components of $G\setminus A$, by Proposition \ref{identity}(a), $\m(\theta,G\setminus A)$ is at least $\m(\theta,D)$, a contradiction.
\end{proof}

The result just proved gives more structure to the Gallai-Edmonds decomposition. For example, $\m(\theta,G)$ is the difference of $|A(G)|$ from the number of components induced by $D(G)$. For the case $\theta=0$, the powerful tools offered by Theorems \ref{main} and \ref{main2} are known as the Gallai-Edmonds Structure Theorem.

\end{document}